\input amstex
\documentstyle{amsppt}
%
\catcode`@=11
\redefine\output@{%
  \def\break{\penalty-\@M}\let\par\endgraf
  \ifodd\pageno\global\hoffset=105pt\else\global\hoffset=8pt\fi  
  \shipout\vbox{%
    \ifplain@
      \let\makeheadline\relax \let\makefootline\relax
    \else
      \iffirstpage@ \global\firstpage@false
        \let\rightheadline\frheadline
        \let\leftheadline\flheadline
      \else
        \ifrunheads@ 
        \else \let\makeheadline\relax
        \fi
      \fi
    \fi
    \makeheadline \pagebody \makefootline}%
  \advancepageno \ifnum\outputpenalty>-\@MM\else\dosupereject\fi
}
\def\Beta{\mathchar"0\hexnumber@\rmfam 42}
\catcode`\@=\active
\nopagenumbers
\chardef\textvolna='176

\chardef\bigalpha='013
\def\negskp{\hskip -2pt}

\def\Ker{\operatorname{Ker}}

\chardef\degree="5E
\def\compos{\,\raise 1pt\hbox{$\sssize\circ$} \,}

\def\blue#1{#1}

\catcode`#=11\def\diez{#}\catcode`#=6
\catcode`&=11\catcode`&=4
\catcode`_=11\def\podcherkivanie{_}\catcode`_=8
\catcode`\^=11\catcode`\^=7
\catcode`~=11\def\volna{~}\catcode`~=\active
\def\mycite#1{\cite{\blue{#1}}\immediate\special{ps:
     ShrHPSdict begin /ShrBORDERthickness 0 def}}
\def\myciterange#1#2#3#4{\cite{\blue{#2#3#4}}\immediate\special{ps:
     ShrHPSdict begin /ShrBORDERthickness 0 def}}
\def\mytag#1{%
    \tag#1}
\def\mythetag#1{\thetag{\blue{#1}}\immediate\special{ps:
     ShrHPSdict begin /ShrBORDERthickness 0 def}}
\def\myrefno#1{\no#1}
\def\myhref#1#2{\blue{#2}\immediate\special{ps:
     ShrHPSdict begin /ShrBORDERthickness 0 def}}
\def\myEarXivlink{\myhref{http://arXiv.org}{http:/\negskp/arXiv.org}}

\def\mytheorem#1{\csname proclaim\endcsname{Theorem #1}}
\def\mytheoremwithtitle#1#2{\csname proclaim\endcsname{Theorem #1#2}}
\def\mythetheorem#1{\blue{#1}\immediate\special{ps:
     ShrHPSdict begin /ShrBORDERthickness 0 def}}
\def\mylemma#1{\csname proclaim\endcsname{Lemma #1}}
\def\mylemmawithtitle#1#2{\csname proclaim\endcsname{Lemma #1#2}}

\def\mycorollary#1{\csname proclaim\endcsname{Corollary #1}}

\def\myconjecture#1{\csname proclaim\endcsname{Conjecture #1}}
\def\myconjecturewithtitle#1#2{\csname proclaim\endcsname{Conjecture #1#2}}

\def\myproblem#1{\csname proclaim\endcsname{Problem #1}}
\def\myproblemwithtitle#1#2{\csname proclaim\endcsname{Problem #1#2}}


\pagewidth{360pt}
\pageheight{606pt}
\topmatter
\title
On some higher degree sign-definite multivariate polynomials 
associated with definite quadratic forms.
\endtitle
\rightheadtext{On some sign-definite multivariate polynomials \dots}
\author
Ruslan Sharipov
\endauthor
\address Bashkir State University, 32 Zaki Validi street, 450074 Ufa, Russia
\endaddress
\email
\myhref{mailto:r-sharipov\@mail.ru}{r-sharipov\@mail.ru}
\endemail
\abstract
     Positive and negative quadratic forms are well known and widely used. 
They are multivariate homogeneous polynomials of degree two taking positive 
or negative values respectively for any values of their arguments not all zero. 
In the present paper a certain higher degree polynomial is associated with 
each quadratic form such that the form is definite if and only if this 
polynomial is sign-definite. 
\endabstract
\subjclassyear{2000}
\subjclass 15A48, 15A45, 11E10\endsubjclass
\endtopmatter
\TagsOnRight
\document

\head
1. Introduction.
\endhead
     Quartic and higher order positive polynomials are of growing interest
(see \myciterange{1}{1}{--}{5}). Trivial examples of them are constructed 
as sums of squares. But in general the polynomial non-negativity is
an NP-hard problem (see \mycite{6}, \mycite{7}). For this reason any 
examples of higher degree sign-definite polynomials are worthwhile.\par
     Let $a(x^1,\ldots,x^n)$ be a quadratic form of $n$ variables\footnotemark.
It is given by the formula
$$
\hskip -2em
a(x^1,\ldots,x^n)=\sum^n_{i=1}\sum^n_{j=1} a_{ij}\,x^i\,x^j. 
\mytag{1.1}
$$
\footnotetext{\,Upper indices for numerating variables in \mythetag{1.1} are 
used according to Einstein's tensorial notation, see \mycite{8}.}%
\adjustfootnotemark{-1}%
The form \mythetag{1.1} is associated with its matrix

$$
\hskip -2em
A=\Vmatrix a_{11}& \hdots & a_{1n}\\
\vdots & \ddots &\vdots\\
a_{n1}& \hdots & a_{nn}\endVmatrix
\mytag{1.2}
$$
which is symmetric, i\.\,e\. $a_{ij}=a_{j\kern 0.2pt i}$. If $x^1,\,\ldots,
\,x^n$ are interpreted as components of a vector in some linear vector space 
$V$, then $a_{ij}$ are components of a twice covariant tensor. Under a linear
change of variables
$$
\hskip -2em
x^i=\sum^n_{j=1}S^i_j\,\tilde x^j,
\mytag{1.3}
$$ 
which is interpreted as a change of basis in $V$, the matrix \mythetag{1.2} is
transformed as 
$$
\hskip -2em
\tilde A=S^{\kern 1pt\sssize\top}\!A\,S. 
\mytag{1.4}
$$
Here $S$ is the transition matrix (see \mycite{8} or \mycite{9}). Its components 
are used in \mythetag{1.3}.\par
      Let $\Lambda$ be a skew-symmetric matrix of the same size as the matrix $A$ 
in \mythetag{1.2}: 
$$
\hskip -2em
\Lambda=\Vmatrix 0 & \lambda_{12}& \hdots & \lambda_{1n}\\
\vspace{1ex}
-\lambda_{12}& 0 & \hdots & \lambda_{2n}\\
\vdots & \vdots &\ddots &\vdots \\
-\lambda_{1n}& -\lambda_{2n} & \hdots & 0\endVmatrix.
\mytag{1.5}
$$
Using the matrices \mythetag{1.2} and \mythetag{1.5}, we define the polynomial
$$
\hskip -2em
P(\lambda_{12},\ldots,\lambda_{n-1\,n})=P(\Lambda)=\det(A-\Lambda). 
\mytag{1.6}
$$
For the definition \mythetag{1.6} to be coordinate covariant we set 
$$
\hskip -2em
\tilde\Lambda=S^{\kern 1pt\sssize\top}\Lambda\,S, 
\mytag{1.7}
$$
which is the same transformation rule as \mythetag{1.4}. From 
\mythetag{1.4} and \mythetag{1.7} we derive 
$$
\hskip -2em
P(\tilde\Lambda)=(\det S)^2\,P(\Lambda). 
\mytag{1.8}
$$\par
     The components $\lambda_{12},\,\ldots,\,\lambda_{n-1\,n}$ of the matrix 
$\Lambda$ in \mythetag{1.6} are interpreted as independent variables, i\.\,e\.
as arguments of the polynomial $P(\Lambda)$. The polynomial $P(\Lambda)$ 
in \mythetag{1.6} is of degree two in each particular variable 
$\lambda_{\kern 0.5pt ij}$. 
However, its total degree is typically higher than two.\par 
     Note that the formula \mythetag{1.6} is somewhat similar to the formula
of the characteristic polynomial of a matrix. Therefore below we shall call 
$P(\Lambda)$ the skew-characteristic polynomial of the form \mythetag{1.1}. 
Studying some properties of this polynomial is the main goal of the present
paper. 
\head
2. Proving the positivity. 
\endhead
\mytheorem{2.1} If a quadratic form with the matrix $A$ is positive,
then its associated skew-characteristic polynomial $P(\Lambda)=\det(A-\Lambda)$
is positive, i\.\,e\. $P(\Lambda)>0$ for any skew-symmetric matrix $\Lambda$. 
\endproclaim
    It is known that the matrix of a positive quadratic form can be brought to 
the unit matrix at the expense of linear transformations of the form 
\mythetag{1.3}, i\.\,e\. in a proper basis (see \mycite{9}). Therefore, relying
on \mythetag{1.8}, without loss of generality we can choose $A=\bold 1$ and 
consider some examples.\par
     {\bf The case} n=2. In this case we easily calculate
$$
\hskip -2em
P(\Lambda)=\vmatrix 1 & -\lambda_{12}\\ 
\vspace{2ex}
\lambda_{12} & 1\endvmatrix=1
+\lambda_{12}^{\kern 2pt 2}>0. 
\mytag{2.1}
$$\par
     {\bf The case} n=3. This case is similar to the previous one:
$$
\hskip -2em
P(\Lambda)=\vmatrix 1 & -\lambda_{12} & -\lambda_{13}\\ 
\vspace{2ex}
\lambda_{12} & 1 & -\lambda_{23}\\
\vspace{2ex}
\lambda_{13} & \lambda_{23} & 1\\
\endvmatrix=1+\lambda_{12}^{\kern 2pt 2}
+\lambda_{13}^{\kern 2pt 2}+\lambda_{23}^{\kern 2pt 2}>0. 
\mytag{2.2}
$$\par
     {\bf The case} n=4. This case is a little bit more complicated than
\mythetag{2.1} and \mythetag{2.2}:
$$
\hskip -2em
\gathered
P(\Lambda)=1+\lambda_{12}^{\kern 2pt 2}+\lambda_{13}^{\kern 2pt 2}
+\lambda_{14}^{\kern 2pt 2}+\lambda_{23}^{\kern 2pt 2}
+\lambda_{24}^{\kern 2pt 2}+\lambda_{34}^{\kern 2pt 2}+\\
+(\lambda_{12}\,\lambda_{34}+\lambda_{23}\,\lambda_{14}
-\lambda_{13}\,\lambda_{24})^{\kern 1pt 2}>0.
\endgathered
\mytag{2.3}
$$
Taking \mythetag{2.1}, \mythetag{2.2}, and \mythetag{2.3} as a background,
we proceed to proving Theorem~\mythetheorem{2.1}.\par
\demo{Proof of Theorem~\mythetheorem{2.1}} Interpreting $x^1,\,\ldots,\,x^n$
in \mythetag{1.1} as the coordinates of a vector $\bold x\in V$, we can 
associate a symmetric bilinear form with the matrix $A$:
$$
\hskip -2em
a(\bold x,\bold y)=\sum^n_{i=1}\sum^n_{j=1} a_{ij}\,x^i\,y^j. 
\mytag{2.4}
$$
The matrix $\Lambda$ in \mythetag{1.5} is not associated with a quadratic form.
However, it is associated with a skew-symmetric bilinear form:
$$
\hskip -2em
\lambda(\bold x,\bold y)=\sum^n_{i=1}\sum^n_{j=1}\lambda_{\kern 0.5pt ij}\,x^i\,y^j. 
\mytag{2.5}
$$
Due to \mythetag{2.4} and \mythetag{2.5} the matrix $A-\Lambda$ in \mythetag{1.6}
is associated with the bilinear form 
$$
\hskip -2em
b(\bold x,\bold y)=a(\bold x,\bold y)-\lambda(\bold x,\bold y),
\mytag{2.6}
$$
which is neither symmetric nor skew-symmetric.\par
     Assume that the matrix $A$ of the positive quadratic form 
$a(x^1,\ldots,x^n)=a(\bold x,\bold x)$ is brought to the unit matrix by choosing
some proper basis in $V$. Then for $\Lambda=0$ in \mythetag{1.6} we have
the following inequality for $P(\Lambda)$: 
$$
\hskip -2em
P(\Lambda)=P(0)=\det(\bold 1)=1>0.
\mytag{2.7}
$$\par
     Further we shall prove that the polynomial $P(\Lambda)$ cannot vanish.
The proof is by contradiction. Indeed, if $P(\Lambda)=0$, then 
$\det(A-\Lambda)=0$ and $A-\Lambda$ is a degenerate matrix. This means that 
the form \mythetag{2.6} has a nonzero kernel\footnotemark:
$$
\hskip -2em
\Ker b=\{\bold x\in V\!:\ b(\bold x,\bold y)=0\ \forall\,y\in V\}\neq\{\bold 0\}.
\mytag{2.8}
$$
\footnotetext{\ Actually the form $b$ has two kernels --- the left kernel and 
the right kernel, both being nonzero. We choose the left kernel in 
\mythetag{2.8} for the sake of certainty.}
\adjustfootnotemark{-1}
Let $\bold x\neq 0$ be a vector belonging to the kernel \mythetag{2.8}. Then
$$
\pagebreak
\hskip -2em
b(\bold x,\bold x)=0. 
\mytag{2.9}
$$
Applying \mythetag{2.6} to \mythetag{2.9} and taking into account that
$\lambda(\bold x,\bold x)=0$ since the bilinear form $\lambda$ 
in \mythetag{2.5} is skew-symmetric, we derive 
$$
a(\bold x,\bold x)=b(\bold x,\bold x)+\lambda(\bold x,\bold x)=0+0=0.
$$
But the equality $a(\bold x,\bold x)=0$ for $\bold x\neq 0$ contradicts
the positivity of the form $a$. The contradiction obtained proves that
$P(\Lambda)$ cannot vanish.\par
     Thus we know that the polynomial $P(\Lambda)$ is a continuous function
of its arguments which is positive for $\Lambda=0$ due to \mythetag{2.7}
and which never vanishes. Therefore $P(\Lambda)$ is always positive. 
Theorem~\mythetheorem{2.1} is proved. 
\qed\enddemo
\head
3. A criterion of definiteness. 
\endhead
\mytheorem{3.1} A quadratic form with the matrix $A$ is definite if and only
if its associated skew-characteristic polynomial $P(\Lambda)=\det(A-\Lambda)$
is sign-definite. 
\endproclaim
\demo{Proof} Assume that the form $a$ is definite. Then it is either positive or
negative. If $a$ is a positive form with the matrix $A$, then $P(\Lambda)>0$,
which follows from Theorem~\mythetheorem{2.1}. If $a$ is negative, then the form
$-a$ is positive. For this form we derive
$$
\hskip -2em
P_{-a}(-\Lambda)=\det(-A+\Lambda)=(-1)^n\,\det(A-\Lambda)=(-1)^n\,P(\Lambda).
\mytag{3.1}
$$
Applying Theorem~\mythetheorem{2.1} to \mythetag{3.1}, we find that $P(\Lambda)>0$
if the dimension $n=\dim V$ is even and $P(\Lambda)<0$ if $n$ is odd. In both cases
$P(\Lambda)$ is sign-definite. This means that the necessity is proved.\par 
     The proof of the sufficiency is by contradiction. Assume that $P(\Lambda)$ 
is sign-definite but the form $a$ is not definite. Then $a$ is either degenerate 
or non-degenerate. If $a$ is degenerate, then $\det A=0$. Choosing $\Lambda=0$ in \mythetag{1.6}, we get $P(0)=\det(A)=0$. The equality $P(0)=0$ contradicts both 
$P(\Lambda)>0$ and $P(\Lambda)<0$.\par
      If $a$ is non-degenerate and indefinite, then its signature is $(m,n-m)$,
where $m\neq 0$ and $n-m\neq 0$. In this case by mean of some proper choice of 
basis in $V$ we can bring the matrix $A$ to the following diagonal form: 
$$
\hskip -2em
A=\Vmatrix 1 & \hdots & 0 & 0 & \hdots & 0\\
\vdots & \ddots & \vdots & \vdots && \vdots\\
0 & \hdots & 1 & 0 & \hdots & 0\\
\vspace{1ex}
0 & \hdots & 0 & -1 & \hdots & 0\\
\vdots && \vdots & \vdots & \ddots & \vdots\\
0 & \hdots & 0 & 0 &\hdots & -1
\endVmatrix.
\mytag{3.2}
$$
Relying on \mythetag{3.2}, we choose the following skew-symmetric matrix 
$\Lambda$:
$$
\hskip -2em
\Lambda=\Vmatrix 0 & \hdots & 0 & 0 & \hdots & 0\\
\vdots & \ddots & \vdots & \vdots && \vdots\\
0 & \hdots & 0 & \lambda_{m\,m+1} & \hdots & 0\\
\vspace{2ex}
0 & \hdots & -\lambda_{m\,m+1} & 0 & \hdots & 0\\
\vdots && \vdots & \vdots & \ddots & \vdots\\
0 & \hdots & 0 & 0 &\hdots & 0
\endVmatrix.
\mytag{3.3}
$$
Substituting \mythetag{3.2} and \mythetag{3.3} into \mythetag{1.6} we
derive
$$
\hskip -2em
P(\Lambda)=(-1)^{n-m-1}\,(-1+\lambda_{m\,m+1}^{\kern 2pt 2}). 
\mytag{3.4}
$$
It is easy to see that the polynomial \mythetag{3.4} is not sign-definite, which 
is again a contradiction. Thus, Theorem~\mythetheorem{3.1} is proved. 
\qed\enddemo
     Theorem~\mythetheorem{3.1} is the main result of the present paper. 
It can be further used as a background in deriving definiteness criteria 
for quartic and higher order forms.\par

\enddocument
\end